\documentstyle[12pt]{article}
\def\Bbb#1{{\mathchoice{\mbox{\bf #1}}{\mbox{\bf #1}}%
{\mbox{$\scriptstyle \bf #1$}}{\mbox{$\scriptscriptstyle \bf #1$}}}}
\def\N{\Bbb N}
\def\R{\Bbb R}
\def\C{\Bbb C}

\def\T{\Bbb T}
\def\Q{\Bbb Q}

\def\cL{{\cal L}}

\begin{document}
\title{Analytic measures and Bochner measurability}
\author{Nakhl\'e H. \ Asmar\ and\ Stephen J.\ Montgomery--Smith\\
Department of Mathematics\\
University of
Missouri--Columbia\\
Columbia, Missouri 65211  U.\ S.\ A.
}
\date{}
\maketitle
			     %
			     %
			     %
			     %
			     %
			     %
			     %
			     %


\section{Introduction}
\newtheorem{weakmble}{Definition}[section]
\newtheorem{weakanalytic}[weakmble]{Definition}
\newtheorem{hypA}[weakmble]{Definition}
\newtheorem{ex1}[weakmble]{Example}
\newtheorem{ex2}[weakmble]{Example}
\newtheorem{ex3}[weakmble]{Example}
\newtheorem{prop5.2}[weakmble]{Proposition}
\newtheorem{cor5.3}[weakmble]{Corollary}
\newtheorem{prop5.5}[weakmble]{Proposition}
\newtheorem{cor5.7}[weakmble]{Proposition}
  Many authors have made
 great strides in extending the celebrated 
F.\ and M.\ Riesz Theorem 
to various abstract settings. 
 Most notably, we have, in chronological order,
the work of Bochner \cite{bochner}, 
Helson and Lowdenslager \cite{hl1}, 
de Leeuw and Glicksberg \cite{dg}, and Forelli
\cite{forelli}.  
These formidable papers build on each other's ideas
and provide broader extensions of the F.\ and M.\ Riesz Theorem.
Our goal in this paper is to  use the
analytic Radon-Nikod\'ym property and prove a representation 
theorem (Main Lemma \ref{mainlemma} below) 
for a certain class of
measure-valued mappings on the real line.
Applications of this result yield the  
the main theorems from \cite{dg} 
and \cite{forelli}.  First, we will review briefly the results
with which we are concerned, and describe our main
theorem.

The F.\ and M.\ Riesz Theorem states that if
a complex Borel measure $\mu$ on the circle is such that
$$\int_{-\pi}^\pi e^{- in t} d\mu(t)=0,\ {\rm for\ all}\ n<0,$$
(i.e.\ $\mu$ is analytic)
then $\mu$ is absolutely continuous with respect 
to Lebesgue measure.  The first extension  
is due to Bochner \cite{bochner} who 
used very elaborate methods to show
that if the Fourier transform of a complex Borel measure 
on the two dimensional
torus vanishes off a sector with opening strictly less than
 $\pi$, then the measure is absolutely continuous
 with respect to the two-dimensional Lebesgue measure.
A few years later, Helson and 
Lowdenslager \cite{hl1}, and de Leeuw and Glicksberg 
\cite{dg} revisited this theorem and offered
different proofs based on their abstract versions
of the F.\ and M.\ Riesz theorem.  
The paper \cite{hl1} is
classical; it contains seminal work in 
harmonic analysis on ordered groups, an  area
of analysis that flourished in the decades that
followed. 
In \cite{hl1}, a measure is called
analytic if its Fourier transform vanishes on the negative
characters, and their version of the F.\ and M.\
Riesz Theorem states:

{\em if a measure $\mu$ is analytic, then its absolutely continuous 
part and its singular part, with respect to 
Haar measure on the group, are both analytic.}

Looking at the F.\ and M.\ Riesz Theorem from a different 
perspective, de Leeuw and Glicksberg 
considered the setting of a compact abelian
group $G$ on which the real line $\R$ is acting by translation via 
a continuous homomorphism from $\R$ into $G$.
Thus the dual homomorphism maps the dual group of $G$ to $\R$.
In this setting, analytic measures are those with
Fourier transforms supported on the inverse image of
the positive real line.
The de Leeuw-Glicksberg version of the F.\ and M.\ Riesz Theorem states:

{\em the Borel subsets of $G$ on which an analytic measure 
vanishes identically is invariant under the action of $\R$.}

De Leeuw and Glicksberg called a measure
whose null sets are invariant under $\R$ quasi-invariant.  With
this terminology, their result states
that every analytic measure is quasi-invariant.

The notion of quasi-invariance and analyticity were extended by Forelli
\cite{forelli} to the setting in which the real line is acting
on a locally compact topological space.
Since Forelli's setting is closest to ours, we 
will describe it in greater detail.

\bigskip
\noindent
{\bf Forelli's main results}  Let $\Omega$ be a locally compact 
Hausdorff space, and let $T:\ t\mapsto T_t$ 
denote a representation
of the real line $\R$ by homeomorphisms of the 
topological space $\Omega$ such that 
the mapping $(t,\omega)\mapsto T_t\omega$ 
is jointly continuous.
The action of $\R$ on $\Omega$ induces,
in a natural way, an action on the Baire measures on
$\Omega$.  With a slight abuse of notation,
if $\mu$ is a Baire measure and $A$ is a Baire
subset of $\Omega$, we write
$T_t\mu$ for the Baire measure whose value at $A$ is
$\mu(T_t A)$.  Denote the Baire subsets by 
$\Sigma$, and the Baire measures 
by $M(\Omega,\Sigma)$, or simply $M(\Sigma)$.
A measure $\nu$ in $M(\Sigma)$ is called quasi-invariant if
the collection of subsets of $\Sigma$ on which 
$\nu$ vanishes identically is invariant by $T$.
That is, $\nu$ is quasi-invariant if
$|\nu|(T_t A)=0 $ for all $t$ if and only if
$|\nu|(A)=0 $.

Using the representation $T$, one can 
define the spectrum of a measure in $M(\Sigma)$ 
(see (\ref{specsbt}) below), which plays 
the role of the support of the Fourier transform of a measure. 
A measure in $M(\Sigma)$ is then called 
analytic if its spectrum lies on the nonnegative real axis.
With this terminology, Forelli's 
main result states that:

{\em  an analytic measure is quasi-invariant.}

As a corollary of this result, 
Forelli \cite[Theorem 4]{forelli} showed that
analytic measures translate continuously.  That is,

{\em  if $\mu$ is analytic, then $t\mapsto T_t\mu$
is continuous from $\R$ into $M(\Sigma)$. }\\

When $\Omega$ is the real line, and $T_t$
stands for translation by $t$, a quasi-invariant measure, or a 
measure for which the mapping
$t\mapsto T_t\mu$ is continuous is necessarily
absolutely continuous with respect to Lebesgue measure.
These facts were observed
by de Leeuw and Glicksberg \cite{dg} and for these reasons
the main results in \cite{dg} and \cite{forelli}
are viewed as extensions of the F.\ and M.\ Riesz Theorem.


\bigskip
\noindent
{\bf Goals of this paper}  
Although Forelli proves that analytic measures 
translate continuously as a consequence of quasi-invariance,
it can be shown that, vice-versa, in the setting of 
Forelli's paper, the quasi-invariance of analytic measures is
a consequence of the continuity of the mapping
$t\mapsto T_t\mu$ (see Section 5 below).
The latter approach is the one that we take in this paper.
As we now describe, this approach
has many advantages, and the main 
results of this paper cannot
be obtained using Forelli's methods.

Let $\Sigma$ denote a 
sigma algebra of subsets of a set $\Omega$ 
and let $M(\Sigma)$
denote the space of complex measures defined on
$\Sigma$.  Suppose that $T:  \ t\mapsto T_t$
is a uniformly bounded group of
isomorphisms of $M(\Sigma)$.  
Using the representation $T$, we can define
the notion of analytic measures as in 
\cite{forelli}, or as described in Definition
\ref{weakanalytic} below.  For an analytic 
$\mu$ in $M(\Sigma)$, we ask:  
under what conditions
on $T$ is the mapping
$t\mapsto T_t\mu$ continuous?
  
Clearly, if this mapping is to be 
continuous, then the
following must hold:
if $\nu$ is analytic such that
for every $A\in \Sigma$,  
$T_t\nu(A)=0$ for almost all $t\in \R$ then $\nu$ is the zero measure.

Our main results (Theorems \ref{riesz1} 
and \ref{minimaintheorem} below)
prove that the converse is also true.
We call the property that we just 
described hypothesis $(A)$ (see 
Definition \ref{hypA}
below), and show, for example, that if 
a representation $T$, given by mappings of the
sigma algebra, satisfies 
hypothesis $(A)$, then
the mapping $t\mapsto T_t\mu$
is Bochner measurable whenever $\mu$ is analytic.   
Using this fact, we can derive with
ease all the main properties of
analytic measures that were obtained by Forelli \cite{forelli}.  
By imposing the right conditions on
$T$, we are able to use the analytic Radon-Nikod\'ym property
of the Banach space $M(\Sigma)$ to give 
 short and perspicuous proofs which dispense
 with several unnecessary conditions on the representation.  
 In particular, in many interesting situations,
 we do not even need the fact that the collection
 of operators $(T_t)_{t\in \R}$ forms a group
 under composition.


\bigskip
\noindent
{\bf Notation and Definitions} \ \ We 
use the symbols $\Q, \R$, and $\C$ 
to denote the rational 
numbers, the real numbers, and the complex numbers respectively.  
The circle group will be denoted by $\T$ and will be 
customarily parametrized as $\{ e^{i t}:\ 0\leq t<2 \pi\}$.
Our measure theory is borrowed 
from Hewitt and Ross \cite{hr1}.  In particular, 
the convolution of measures and functions is defined as in
\cite[\S 20]{hr1}.  We denote by
$M(\R)$
the Banach space of complex regular Borel measures on $\R$.
The space of  
Lebesgue measurable integrable functions on $\R$ is denoted
by $L^1(\R)$, and the space of essentially 
bounded measurable functions by $L^\infty(\R)$. 
The spaces $H^1(\R)$ and $H^\infty(\R)$ are defined as follows:
$$H^1(\R)=\left\{
f\in L^1(\R):\ \widehat{f}(s)=0,\ s\leq 0
\right\};$$
and
$$H^\infty(\R)=\left\{
f\in L^\infty(\R): \ \int_\R f(t) g(t) dt=0\ {\rm for\ all\ }
g\in H^1(\R)
\right\}.$$
Let $(\Omega, \Sigma)$ denote a 
measurable space and let $\cL^\infty(\Sigma)$ denote the 
bounded measurable functions on $\Omega$.  Denote 
by $M(\Sigma)$ the Banach space of countably additive 
complex measures on $(\Omega, \Sigma)$ with the 
total variation norm.  
Suppose that $T=(T_t)_{t\in\R}$ is a collection
of uniformly bounded 
invertible isomorphisms of $M(\Sigma)$ with
\begin{equation}
\|T_t^{\pm 1}\|\leq c
\label{uniformlybded}
\end{equation}
for all $t\in \R$, where $c$ is a positive constant.  
(Note that we do not require that $(T_t)^{-1}=T_{-t}$, but only
that $T_t$ be invertible.)
The following definition determines the class of measures 
that we will be studying.
\begin{weakmble}
Let $(T_t)_{t\in \R}$ be as above. 
A measure $\mu\in M(\Sigma)$ is 
called weakly measurable if 
for every $A\in \Sigma$, the mapping $t\mapsto T_t\mu(A)$ is 
Lebesgue measurable on $\R$.
\label{weakmble}
\end{weakmble}
We next introduce our notion of analyticity.  We will show at the 
end of this section that our notion of analyticity agrees
with Forelli's notion in \cite{forelli},
when restricted to Forelli's setting.  
\begin{weakanalytic}
Let $(T_t)_{t\in\R}$ be a  
uniformly bounded collection of 
isomorphisms of $M(\Sigma)$.  
A weakly measurable $\mu\in M(\Sigma)$ is 
called weakly analytic if the mapping $t\mapsto T_t\mu(A)$ 
is in $H^\infty(\R)$ for every $A\in\Sigma$.\\
\label{weakanalytic}
\end{weakanalytic}
Our Main Theorem (Theorem \ref{minimaintheorem} below) 
states that, under a certain condition on $T$
that we described in our introduction, if
$\mu$ is weakly analytic then the mapping
$t\mapsto T_t\mu$ is Bochner measurable.  This key property
is presented in the following definition.
\begin{hypA}
Let $T=(T_t)_{t\in\R}$ be a  
uniformly bounded collection of 
isomorphisms of $M(\Sigma)$. 
Then $T$ is said   
to satisfy hypothesis $(A)$ if whenever $\mu$ is weakly analytic
in $M(\Sigma)$, such that for every $A\in \Sigma$,  
$T_t\mu(A)=0$ for almost all
$t\in \R$, then $\mu$ must be the zero measure.
\label{hypA}
\end{hypA}
We emphasize here that the set of $t's$ for which the equality
$T_t\mu(A)=0$ holds depends in general on $A$.  
Hypothesis $(A)$ is crucial to our
study.  
We offer two main sources of examples where it is satisfied.  
The first one 
is related to Forelli's
setting \cite{forelli}.
\begin{ex1}
{\rm Suppose that $\Omega$ is a topological 
space and $\left( T_t\right)_{t\in\R}$ is a collection 
of homeomorphisms of $\Omega$ onto itself 
such that the mapping
$$(t,\omega)\mapsto T_t\omega$$
is jointly continuous.  Let $\Sigma$ denote the 
Baire subsets of $\Omega$
(hence $\Sigma$ is the smallest 
$\sigma$-algebra such that all   
the continuous complex-valued functions 
are measurable with respect to $\Sigma$.)  This 
is Forelli's setting, except that we 
do not require from $\Omega$ to be a locally compact 
Hausdorff space, and more interestingly, we 
do not assume (thus far) that $(T_t)_{t\in\R}$ forms a group.  
For any Baire measure $\mu$, define $T_t\mu$ on the Baire
sets by $T_t\mu(A)=\mu(T_t(A))$.  
Now, suppose that
$\mu$ is such that $T_t\mu(A)=0$ for almost all
$t$, for any given Baire set $A$.  Then it follows that for any
bounded continuous function $f$\ that
$\int f\circ T_t \,d\mu = 0$\ for almost all $t$.  Since
the map $t\mapsto \int f\circ T_t \,d\mu$\ is continuous, it follows
that $\int f\,d\mu = 0$.  Now suppose that $A = f^{-1}(0,\infty)$.  Then
$\mu(A) = \lim_{n\to\infty} \int\max\{0,\min\{f,1\}\}^{1/n} \, d\mu = 0$.
From this, it is easy to conclude that $\mu = 0$, and so
$T$ satisfies hypothesis $(A)$}.
\label{ex1}
\end{ex1}
Our second source of examples is given by
the abstract Lebesgue spaces which provide 
ideal settings to study analytic measures, in the sense that
the main results of this paper hold with very relaxed conditions on the 
representation.  (See Theorem \ref{riesz1} and Remarks
\ref{8nov94remark} and \ref{8nov94remark2} below.)
\begin{ex2}
{\rm Suppose that $\Sigma$ is countably generated.  Then 
any uniformly bounded collection $(T_t)_{t\in\R}$  
by isomorphisms of $M(\Sigma)$ satisfies hypothesis $(A)$.  
The proof follows easily from definitions.
}
\label{ex2}
\end{ex2}
The next example will be used to 
construct counterexamples 
when a representation 
fails hypothesis $(A)$.  It also serves to 
illustrate the use of hypothesis $(A)$.
\begin{ex3}
{\rm
(a)  Let $\Sigma$ denote the sigma algebra of 
countable and co-countable subsets of $\R$.  
Define $\nu\in M(\Sigma)$ by
$$
\nu (A)=\left\{
\begin{array}{ll}
1 & \mbox{if $A$ is co-countable,}\\
0 &  \mbox{if $A$ is countable.}
\end{array}
\right.
$$
Let $\delta_t$ denote the point mass at $t\in \R$; 
and take $\mu=\nu-\delta_0$.  Consider the 
representation $T$ of $\R$ given by translation by $t$.  
Then it is easily verified that 
$\|\mu\|>0$, whereas for every $A\in \Sigma$\ we have that 
$T_t(\mu)(A)=0$ 
for almost all $t\in \R$.
Hence the 
representation $T$ does not satisfy hypothesis $(A)$.

The following generalization of (a) will be 
needed in the sequel.\\
(b)  Let $\alpha$ be a real number and let 
$\Sigma$, $\mu$, $\nu$, $\delta_t$, and $T_t$ 
have the same meanings as in (a).  Define a 
representation $T^\alpha$ by
$$T_t^\alpha = e^{i\alpha t}T_t.$$ 
Arguing as in (a), it is easy to see that $T^\alpha$ 
does not satisfy hypothesis $(A)$.
}
\label{ex3}
\end{ex3}

\bigskip
\noindent
{\bf Organization of the paper} 
In the rest of this section, we introduce some
notions from spectral synthesis of bounded functions and
show how our definition of analytic measures compares to
Forelli's notion.  
Section 2 contains our Main Lemma and some preliminary applications
to generalized analyticity.
Although this section does not contain our most general results,
it shows the features of our new approach which 
is based on 
the analytic Radon-Nikod\' ym property of 
Bukhvalov and Danilevich \cite{bd}.  
In Section 3, we deal with a one-parameter
group acting on $M(\Sigma)$. 
Using results from Section 2, we derive 
our main application which concerns the Bochner measurability
of the mapping $t\mapsto T_t\mu$.  
In Section 4, we specialize our study 
to representations that are defined by mappings of the sigma 
algebra and prove results concerning the 
Lebesgue decomposition of analytic measures. 
Finally in Section 5, we assume that the representation
is given by point mappings 
and give a short and simple proof
that analytic measures are
quasi-invariant.  The results of Sections 4 and 5 generalize their
counterparts in 
Forelli's paper.  
We also show by examples that Forelli's
approach cannot possibly imply the results of the earlier sections.
Section 5 concludes with remarks about further extensions of
our methods to the setting where $\R$ is replaced by any
locally compact abelian group with an ordered dual group.
These extensions combine the version of the F. and M. Riesz due to 
Helson and Lowdenslager \cite{hl1} with the results of this paper.

\bigskip

Now let us discuss the definition of analyticity
according to Forelli.  We give this definition
in our general setting of a representation $T$ of $\R$ acting
on $M(\Sigma)$.  For a weakly measurable $\mu\in M(\Sigma)$, we let 
$${\cal J}(\mu)=
  \left\{f\in L^1(\R):\ \ \int_\R T_t \mu(A) f(s-t) dt =0
  \mbox{\ for almost all $s\in\R$\ for all $A\in\Sigma$}\right\}.$$
Define the $T$-spectrum of $\mu$ by
\begin{equation}
{\rm spec}_T (\mu)= \bigcap_{f\in {\cal J}(\mu)}
\left\{
\chi\in\R:\ \ \widehat{f}(\chi)=0
\right\}.
\label{specsbt}
\end{equation}
The measure $\mu$ is called $T$-analytic if its
$T$-spectrum is contained in $[0,\infty)$.

The result we need to equate Forelli's notion of analyticity with the notion
we present is the following.
\begin{cor5.7}
A measure $\mu\in M(\Sigma)$ is weakly analytic if and only if 
it is $T$-analytic.
\label{cor5.7}
\end{cor5.7}
It follows almost immediately from the definitions that if $\mu$\ is
weakly analytic, then it is $T$-analytic.  The converse is not so
obvious, and requires the following notions.  
Our reference for the rest of this section is Rudin
\cite[Chapter 7]{rudin}.

Given $\phi\in L^\infty(\R)$, define its ideal by
$${\cal J}(\phi)=
  \left\{f\in L^1(\R):\ \ f * \phi = 0\right\}.$$
One definition of the spectrum of $\phi$ is (see \cite[Chapter 7, Theorem
7.8.2]{rudin})
\begin{equation}
\sigma(\phi) = \bigcap_{f\in {\cal J}(\phi)}
\left\{
\chi\in\R:\ \ \widehat{f}(\chi)=0
\right\}.
\label{spec}
\end{equation}
A set $S \subset \R$\ is called a set of spectral synthesis if whenever
$\phi \in L^\infty(\R)$\ with $\sigma(\phi) \subset S$, then $\phi$\
can be 
approximated in 
the weak-* topology of $L^\infty(\R)$ by linear 
combinations of characters from $S$.  
(See \cite[Section 7.8]{rudin}.)  With this definition,
the following proposition follows easily.
\begin{prop5.2}
Let $S$ be a nonvoid closed subset of $\R$ that is a set of spectral
synthesis.  If $\phi\in L^\infty(\R)$\ with $\sigma(\phi) \subset S$,
then
\[
\int_\R f(x)g(x)dx=0
\]
for all $g$ in $L^1(\R)$ such that $\widehat{g}= 0$ on $-S$.
\label{prop5.2}
\end{prop5.2}
The proof of Proposition \ref{cor5.7} now follows immediately from the
fact that $[0,\infty)$\ is a set of spectral synthesis 
(see \cite[Theorem 7.5.6]{rudin}).

\section{The Main Lemma}
\newtheorem{defarnp}{Definition}[section]
\newtheorem{essbound}[defarnp]{Lemma}
\newtheorem{mainlemma}[defarnp]{Main Lemma}
\newtheorem{cormainlemma}[defarnp]{Theorem}
\newtheorem{riesz1}[defarnp]{Theorem}
\newtheorem{propertyofg}[defarnp]{Theorem}
In our proofs we use the notions of Bochner measurability
and Bochner integrability.  A function $f$ from a $\sigma$-finite
measure space $(\Omega,\Sigma,\mu)$
to a Banach space $X$ is Bochner measurable if it satisfies one of the
following two, equivalent, conditions:
\begin{itemize}
\item
$f^{-1}(A)\in\Sigma$ for any open subset $A$ of $X$, and there is a
set $E \in \Sigma$\ such that $\mu(\Omega\setminus E) = 0$\ and $f(E)$\
is separable;
\item
there is a sequence of simple functions $f_n:\Omega\to X$\ such that
$f_n \to f$\ a.e.
\end{itemize}
Furthermore, if $\int \| f \| d\mu < \infty$, then we say that $f$
is Bochner integrable, and it is possible
to make sense of $\int f d\mu $ as an element of $X$.  In particular,
if $P:X\to Y$ is a bounded 
linear operator between two Banach spaces, then
$P\left(\int f d\mu\right) = \int Pf d\mu$.
We refer the reader to \cite[Section 3.5]{hp}.

In this section we prove our Main Lemma about
the Bochner measurability of functions
defined on $\R$ with values in a Banach space with the
analytic Radon-Nikod\'ym property.
This property of Banach spaces was introduced by
Bukhvalov (see, for example \cite{bukh}) to extend
the basic properties of functions in the Hardy spaces on the disc
to vector-valued functions.  

Let ${\cal B}$ denote the Borel subsets of $\T$, and let $X$ denote 
a complex Banach space.  A vector-valued measure 
$\mu:\  {\cal B}\rightarrow X$ of bounded variation (in symbols, 
$\mu\in M({\cal B},X))$ is called 
analytic if 
$$\int_0^{2\pi} e^{-i n t} d\mu (t) = 0\ \ {\rm for\ all}\ n < 0.$$
Analytic measures were extensively studied by Bukhvalov 
and Danilevich (see for example \cite{bd}).  We owe to them the following
definition.
\begin{defarnp}
A complex Banach space $X$ is said to have the analytic 
Radon-Nikod\'ym property (ARNP) if every $X$-valued analytic 
measure $\mu$ in $M({\cal B},X)$ has a Radon-Nikod\'ym 
derivative ---  that is, there is a Bochner measurable $X$-valued 
function $f$ 
in the space of Bochner integrable functions, $L^1(\T,X)$, 
such that
$$\mu(A)=\int_A f dt$$
for all $A\in {\cal B}$.
\label{defarnp}
\end{defarnp}
Like the ordinary 
Radon-Nikod\'ym property, the analytic
Radon-Nikod\'ym property is
about the existence of a Bochner measurable derivative for
vector-valued measures.  However, the difference between the two 
properties, due to the fact that ARNP concerns only analytic measures,
makes the class of Banach spaces with the ARNP strictly larger 
than the class of Banach
spaces with the Radon-Nikod\'ym property.
In this paper, all we need from this theory is the
basic fact that $M(\Sigma)$ has ARNP.  Here, as before, $M(\Sigma)$
denotes the Banach space of complex measures on an arbitrary
$\sigma$-algebra $\Sigma$ of subsets of a set $\Omega$.
According to 
\cite[Theorem 1]{bd},
a Banach lattice $X$ has ARNP if and only if $c_0$ does not embed
in $X$.  (Here, as usual, $c_0$ denotes the linear space of 
complex sequences tending to zero at infinity, and Banach
lattices can be real or complex.)
Since $M(\Sigma)$ is a Banach lattice 
that does not contain a copy of $c_0$, it follows that 
$M(\Sigma)$ has the analytic Radon-Nikod\'ym property.  
(To see that $c_0$ does not embed in $M(\Sigma)$, 
note that $M(\Sigma)$ is weakly sequentially complete, but that
$c_0$ is not.  See \cite[Theorem IV.9.4]{df}.)

Before we state our lemma, we
describe the setting in which it will be used. This will
clarify its statement and proof.

Let $E$ denote the subspace of $\cL^\infty(\Sigma)$ 
consisting of the bounded simple functions on $\Omega$.  
The subspace $E$ embeds isometrically in $M(\Sigma)^*$.  
Then $E$ is a norming subspace of $M(\Sigma)^*$\ for $M(\Sigma)$.
It is also easy to verify that every weak-* sequentially 
continuous linear functional on $E$ is given by point 
evaluation.  That is, if $L:\ E\rightarrow \C$ is 
weak-* sequentially continuous, then there is a 
measure $\mu\in M(\Sigma)$ such that 
\begin{equation}
L(\alpha)=\int_\Omega \alpha \,d\mu
\label{lofe}
\end{equation}
for every $\alpha\in E$.  To verify this fact, it is enough to 
show that the set function given by
\begin{equation}
\mu(A)=L(1_A)
\label{mua}
\end{equation}
defines a measure in $M(\Sigma)$, and this is 
a simple consequence of the
weak-* sequential continuity of $L$.
  
Let $T=(T_t)_{t\in\R}$ be a family of 
uniformly bounded isomorphisms of $M(\Sigma)$
such that (\ref{uniformlybded}) holds.  
Suppose that $\mu\in M(\Sigma)$ is weakly analytic,
and let $f(t)=T_t\mu$ for all $t\in\R$.
Then $\|f(t)\|\leq c\|\mu\|$ 
where $c$ is as in (\ref{uniformlybded}), and
for all $\alpha\in E$, the function $t\mapsto \alpha(f(t))$
is in $H^\infty(\R)$.  With this setting in mind, we
state and prove our Main Lemma.
\begin{mainlemma}
Suppose that $X$ is a complex Banach space with the 
analytic Radon-Nikod\'ym property, and that $E$ is a norming 
subspace of $X^*$, the Banach dual space of $X$.  
Suppose that for every weak-* sequentially continuous functional 
\begin{equation}
L:\ E\rightarrow \C
\label{functional}
\end{equation}
there is an element $x\in X$ such that 
\begin{equation}
L(\alpha)=\alpha(x)
\label{lofeinlemma}
\end{equation}
for all $\alpha\in E$.  Let $f:\ \R\rightarrow X$ be such that 
\begin{equation}
\sup_t\|f(t)\|<\infty
\label{fbded}
\end{equation}
and
\begin{equation}
t\mapsto \alpha(f(t))
\label{fmble}
\end{equation}
is a Lebesgue measurable function in $H^\infty(\R)$ 
for all $\alpha\in E$.  Then there is a 
Bochner measurable function, essentially bounded
\begin{equation}
g:\ \R\rightarrow X
\label{gft}
\end{equation}
 such that for every $\alpha\in E$, we have 
\begin{equation}
\alpha(g(t))=\alpha(f(t))
\label{mainequality}
\end{equation}
for almost all $t\in \R$.  (The set of $t$'s 
for which (\ref{mainequality}) 
holds may depend on $\alpha$.)
\label{mainlemma}
\end{mainlemma}
\begin{essbound}
Suppose that $X$\ is a Banach space, and that $G:\T \to X$\ is a
Bochner integrable function for which there is a constant
$c$\ such that for all Borel sets $A\subset \T$\
\begin{equation}
\left\| \int_A G(\theta) {d\theta\over 2\pi} \right\| \le c \lambda(A) .
\label{essboundeqn}
\end{equation}
(Here $\lambda(A)$\ denotes the Lebesgue measure of $A$).
Then $G$\ is essentially bounded.
\label{essbound}
\end{essbound}
{\bf Proof.}  There is a function $H:\T\to X$\ such that $H = G$\ a.e., and
the range of $H$\ is separable.  Thus there is a countable sequence
$\{\alpha_n\}\subset X^*$\ such that $\|\alpha_n\| \le 1$\ and $\|H(\theta)\| = 
\sup_n \alpha_n(H(\theta))$\ for all $\theta \in \T$.
From (\ref{essboundeqn}), it immediately follows that for
every $n \in \N$\ that
$$ \int_A \alpha_n(H(\theta)) {d\theta\over 2\pi} \le c \lambda(A) ,$$
from whence it follows that $\alpha_n(H(\theta)) \le c$\ a.e.  Hence
$\|H(\theta)\| = \sup_n \alpha_n(H(\theta)) \le c$\ a.e., 
and the result follows.

\bigskip
\noindent
{\bf Proof of Main Lemma \ref{mainlemma}.}  Let $\phi (z)=i \frac{1 - z}{1+z}$
be the conformal mapping of the unit disk onto 
the upper half plane, mapping $\T$\ onto $\R$.  
Let $F=f\circ \phi$.  For every $\alpha\in E$ we have that
$\theta\mapsto\alpha (F(\theta))\in H^\infty (\T)$, since by assumption
$\alpha (f(t))\in H^\infty (\R)$.
Consequently, we have
\begin{equation}
\int_0^{2 \pi} \alpha (F(\theta))e^{i n \theta} 
\frac{d \theta}{2 \pi}= 0 \ \
{\rm for\ all}\ n > 0.
\label{negativefc}
\end{equation}
For $\alpha\in E$, define a measure $\mu_\alpha$ on 
the Borel subsets of $\T$ by
\begin{equation}
\mu_\alpha(A)=\int_A \alpha(F(\theta))\frac{d\theta}{2 \pi}.
\label{musbe}
\end{equation}
Then for all continuous functions $h$ on $\T$, we have
\begin{equation}
\int_\T h(\theta) d\mu_\alpha(\theta)=\int_\T h(\theta)
\alpha(F(\theta))\frac{d\theta}{2\pi}.
\label{intmusbe}
\end{equation}
We 
now claim that for every $A\in {\cal B}(\T)$, 
the mapping $\alpha\mapsto\mu_\alpha(A)$ is weak-* 
sequentially continuous.  That is, if $\alpha_n\rightarrow \alpha$ 
in the weak-* topology of $E$, then
$\mu_{\alpha_n}(A)\rightarrow\mu_\alpha(A)$.  To prove 
this claim, note that if $\alpha_n\rightarrow \alpha$ weak-*, 
then we have that for all $\theta\in \T$, 
\begin{equation}
\alpha_n(F(\theta))\rightarrow \alpha(F(\theta)).
\label{bdedconv}
\end{equation}
Also, by the uniform boundedness principle, we have that
$$\sup_n \|\alpha_n\|=M<\infty .$$ 
So, for all $\theta\in\T$, we have 
$|\alpha_n(F(\theta))|\leq \|\alpha_n\| \|F(\theta)\|
\leq C$.  Hence by bounded convergence, it 
follows from (\ref{bdedconv}) that
\begin{equation}
\mu_{\alpha_n}(A)=
\int_A \alpha_n(F(\theta))\frac{d\theta}{2 \pi}\rightarrow
\int_A \alpha(F(\theta))\frac{d\theta}{2 \pi}
=\mu_\alpha(A),
\label{weak-*conv}
\end{equation}
establishing the desired weak-* sequential continuity.  
By the hypothesis of the lemma, there 
is $\mu (A)\in X$ such that the mapping
$\alpha\mapsto \mu_\alpha(A)$ is given by    
\begin{equation}
\alpha\mapsto \mu_\alpha(A)=\alpha(\mu(A))\ \ {\rm for\ all}\ \alpha\in E.
\label{musbea}
\end{equation}
We now show that for all Borel subsets $A$\ of $\T$\ that
\begin{equation}
\label{muineq}
\|\mu(A)\| \le \int_A \| F(\theta)\| {d\theta\over2\pi} .
\end{equation}
This follows, because $E$\ is norming, and hence, given $A$, and $\epsilon>0$, 
there is an $\alpha\in E$\ with
$\|\alpha\| \le 1$\ and $\|\mu(A)\| \leq |\alpha(\mu(A))| + \epsilon$, and
because
$$ |\alpha(\mu(A))| = \left|\int_A \alpha(F(\theta)) {d\theta\over 2\pi}\right|
   \le \int_A \| F(\theta)\| {d\theta\over2\pi} .$$
From (\ref{muineq}), it is easily seen that the set mapping $A\mapsto \mu(A)$ 
defines an $X$-valued measure of bounded variation on 
the Borel subsets of $\T$.  Let $n$ 
be a positive integer and let $\alpha\in E$.  We have
$$
\alpha\left(\int_\T e^{in\theta}d\mu(\theta)
\right)                 =     \int_\T e^{in\theta} d\mu_\alpha(\theta)
			=     \int_\T e^{in\theta} 
				\alpha(F(\theta))\frac{d\theta}{2 \pi}
			=     0.
$$
And so, since $E$ is norming, it follows that 
$$\int_\T e^{in\theta}d\mu(\theta)=0 $$
for all $n>0$.  Now, appealing to the analytic 
Radon-Nikod\'ym property of $X$, we find a 
Bochner integrable function $G:\ \ \T\rightarrow X$ such that
\begin{equation}
\mu(A)=\int_A G(\theta)\frac{d\theta}{2 \pi}
\label{eq1.24}
\end{equation}
for all Borel subsets $A$ of $\T$.  
Using (\ref{eq1.24}), (\ref{musbea}), and (\ref{musbe}), 
we see that, for all $\alpha\in E$ and all $A\in {\cal B}$,
$$
\int_A \alpha(G(\theta))\frac{d\theta}{2 \pi}        =     \alpha(\mu(A))       
						=     \mu_\alpha(A)        
			=     \int_A \alpha(F(\theta))\frac{d\theta}{2 \pi}.
$$
Since this holds for all $\alpha\in E$ and all $A\in {\cal B}$, 
we conclude that, for a given $\alpha\in E$, 
\begin{equation}
\alpha(G(\theta))=\alpha(F(\theta))\ \ \mbox{a.e.}\ \theta .
\label{egtheta}
\end{equation}
From (\ref{muineq}), (\ref{eq1.24}) and Lemma \ref{essbound}, it follows
that $G$\ is essentially bounded.
Let
\begin{equation}
g(t)=G\left(\phi^{-1}(t)\right).
\label{defofg}
\end{equation}
Then $g$ is Bochner measurable, essentially bounded, and for each
$\alpha\in E$, 
for almost all $t\in \R$,
$$\alpha(g(t))=\alpha(G(\phi^{-1}(t)))=\alpha(F(\phi^{-1}(t)))
=\alpha(f(t)),$$
completing the proof.

\bigskip

When applied in the setting that we described
before the lemma, we obtain the following important consequence.
\begin{cormainlemma}
Let $(T_t)_{t\in\R}$ be a one-parameter family of 
uniformly bounded operators on $M(\Sigma)$ satisfying
(\ref{uniformlybded}), and 
let $\mu$ be a weakly analytic measure in $M(\Sigma)$.  
Then there is a Bochner measurable, essentially bounded function 
$g:\R\rightarrow M(\Sigma)$ such that, for every $A\in \Sigma$, 
$$g(t)(A)=T_t\mu(A),$$
for almost all $t\in \R$.
\label{cormainlemma}
\end{cormainlemma}
Note that the set of $t$'s 
for which the equality in this theorem holds depends on $A$.   
Our goal in the next section is to establish this equality
for all $A\in \Sigma$ and almost all $t\in \R$,
under additional conditions on $T$.
This will imply that the mapping $t\mapsto T_t\mu$
is Bochner measurable 
when $\mu$ is weakly analytic.  
However, when the sigma algebra
is countably generated, this result 
is immediate without any further assumptions on the 
representation. We state it here for ease of reference.
\begin{riesz1}
Suppose that $\Sigma$ is countably generated, and let 
$(T_t)_{t\in\R}$ be a one-parameter family of 
isomorphisms of $M(\Sigma)$ satisfying
(\ref{uniformlybded}).  Suppose that 
$\mu$ is a weakly analytic measure in $M(\Sigma)$.  
Then there is a Bochner measurable function 
$g:\ \R\rightarrow M(\Sigma)$ such that, 
$$g(t)=T_t\mu$$
for almost all $t\in\R$.
\label{riesz1}
\end{riesz1}
{\bf Proof.}  Suppose that $\Sigma=\sigma(\{A_n\}_{n=1}^\infty)$, where
the set $\{A_n\}_{n=1}^\infty$\ is closed under finite unions and
intersections.
Apply Theorem \ref{cormainlemma} to obtain a Bochner 
measurable function $g$ from $\R$ into $M(\Sigma)$ such that
for almost all $t\in \R$ and all $n$ we have
$T_t\mu(A_n)=g(t)(A_n)$.  Since $\Sigma$\ is the closure of 
$\{A_n\}_{n=1}^\infty$ 
under nested countable unions and intersections,
the theorem follows.

\section{Analyticity of measures and Bochner measurability}
\newtheorem{propertyg2}{Proposition}[section]
\newtheorem{firstrelation}[propertyg2]{Lemma}
\newtheorem{secondrelation}[propertyg2]{Lemma}
\newtheorem{minimaintheorem}[propertyg2]{Main Theorem}
\newtheorem{poissonlimit}[propertyg2]{Theorem}
\newtheorem{continuity}[propertyg2]{Theorem}
\newtheorem{stronganalytic}[propertyg2]{Definition}
\newtheorem{spect=spectx}[propertyg2]{Lemma}
\newtheorem{corbigmaintheorem2}[propertyg2]{Theorem}
\newtheorem{8nov94remark}[propertyg2]{Remark}
\newtheorem{applem1}[propertyg2]{Lemma}
\newtheorem{appcor1}[propertyg2]{Corollary}
\newtheorem{appcor2}[propertyg2]{Corollary}

In this section, we prove our main result which states that
if $\mu$ is weakly analytic, then the mapping
$t\mapsto T_t\mu$ is Bochner measurable from $\R$ into
$M(\Sigma)$.  As an immediate consequence of this result
we will obtain that the Poisson integral of a
weakly analytic measure converges in 
$M(\Sigma)$ to the measure, and we also obtain that the 
mapping $t\mapsto T_t\mu$ is continuous.
Both of these results are direct analogues of
classical properties of analytic measures on the real
line.

The proofs in this section require the use of
convolution.  To define this
operation and to derive its basic properties, we will
need additional conditions on the representation $T$.
We start by stating these conditions, setting in the 
process the notation for this section.  These conditions are automatically
satisfied in the case of a representation by mappings of the 
given sigma-algebra.  

We let $T=(T_t)_{t\in\R}$ be a one-parameter group 
of isomorphisms of $M(\Sigma)$\ for which (\ref{uniformlybded}) holds,
satisfying hypothesis $(A)$.  
Here, as before, $M(\Sigma)$ is the Banach space of countably
additive complex measures on an arbitrary sigma algebra $\Sigma$
of subsets of a set $\Omega$.  
We will suppose throughout 
this section that the adjoint of $T_t$ maps $\cL^\infty$ into itself; in
symbols:  
\begin{equation}
T_t^* :\ \cL^\infty (\Sigma) \rightarrow \cL^\infty (\Sigma).
\label{adjointcond}
\end{equation}
Although this property will not appear explicitly in the proofs of the 
main results, we use it at the end of this 
section to establish basic properties of
convolutions of Borel measures on $\R$ with weakly measurable
$\mu$ in $M(\Sigma)$.  
More explicitly, suppose that $\nu\in M(\R)$ and $\mu$ is weakly measurable.
Define a measure $\nu*_T\mu$ (or simply
$\nu *\mu$, when there is no risk of confusion) 
in $M(\Sigma)$ by
$$\nu*_T\mu(A)=\int_\R T_{-t} \mu (A) d\nu (t),$$
for all $A\in \Sigma$.  It is the work of a moment to show that
this indeed defines a measure in $M(\Sigma)$.  We have:
\begin{itemize}
\item for all $t\in \R$, $T_t(\nu*_T\mu)=\nu*_T(T_t\mu)$;
\item the measure $\nu*_T\mu$ is weakly measurable;
\item for $\sigma,\nu\in M(\R)$, and $\mu$ weakly measurable,
$\sigma*_T(\nu*_T\mu)=(\sigma*\nu)*_T\mu$.
\end{itemize}
For clarity's sake, we postpone the proofs of these
results until the end of the section, and proceed towards
the main results.  

Recall the definition of the 
Poisson kernel on $\R$:  for $y>0$, let
$$P_y(x)=\frac{1}{\pi} \frac{y}{x^2+y^2}\ ,$$
for all $x\in \R.$  Let $\mu$ be a weakly analytic measure
in $M(\Sigma)$, and let $g$ be the Bochner
measurable function defined on $\R$ with values in  
$M(\Sigma)$, given by Lemma \ref{mainlemma}.  Form the Poisson
integral of $g$ as follows
\begin{equation}
P_y* g(t)=\int_\R g(t-x) P_y(x) dx,
\label{piofg}
\end{equation}
where the integral exists as a Bochner integral.    
Because the function $g$ is essentially bounded, we have the following result
whose proof follows as in the classical 
setting for scalar-valued functions.
\begin{propertyg2}
With the above notation, we 
have that 
\begin{equation}
\lim_{y\rightarrow 0} P_y*g(t)= g(t)
\label{8nov94(2)}
\end{equation}
for almost all $t\in \R$.
\label{propertyg2}
\end{propertyg2}
We can now establish basic relations between 
the Poisson integral of the function $g$ and 
the measure $\mu$.    
\begin{firstrelation}
For all $t\in\R$, we have
$$P_y*g(t)=P_y*T_t\mu .$$
\label{firstrelation}
\end{firstrelation}
{\bf Proof.} For $A\in\Sigma$, we have
\begin{eqnarray*}
P_y*g(t)(A)     &=&     \int_\R g(s)(A)P_y(t-s)ds\\
		&=&     \int_\R\left(T_s\mu\right)(A)P_y(t-s)ds\\
		&=&     \int_\R\left(T_{t-s}\mu\right)(A)P_y(s)ds\\
		&=&     P_y*\left(T_t\mu\right)(A).
\end{eqnarray*}
Since this holds for all $A\in \Sigma$, the lemma follows.
\begin{secondrelation}
Let $t_0$ be any real number such that
(\ref{8nov94(2)}) holds.  Then, for all $t\in \R$, we have
$$\lim_{y\rightarrow 0} P_y*T_t\mu=T_{t-t_0}\left(g(t_0)\right),$$
in $M(\Sigma)$.
\label{secondrelation}
\end{secondrelation}
{\bf Proof.} 
Since $P_y*g(t_0)\rightarrow g(t_0)$, it follows that
$$T_{t-t_0}\left( P_y*g(t_0)\right) 
\rightarrow T_{t-t_0}\left(g(t_0)\right).$$
Using Lemmas \ref{firstrelation} and the basic
properties of convolutions, we get
$$T_{t-t_0}\left( P_y*g(t_0)\right) = T_{t-t_0}\left(P_y*T_{t_0}\mu\right)
=P_y*T_t\mu,$$
establishing the lemma.

\bigskip

We can now prove the main result of this section.
\begin{minimaintheorem}
Suppose that $T=(T_t)_{t\in\R}$ is a
 group of isomorphisms of $M(\Sigma)$ 
satisfying hypothesis $(A)$ and such that
(\ref{uniformlybded}) and (\ref{adjointcond}) hold.
Let $\mu$ be a weakly analytic measure, and let $g$ 
be the Bochner measurable function on $\R$ 
constructed from $\mu$ as in Theorem \ref{cormainlemma}.  Then
for almost all $t\in\R$, we have
$$T_t\mu=g(t).$$
Consequently, the mapping $t\mapsto T_t\mu$ is Bochner measurable.
\label{minimaintheorem}
\end{minimaintheorem}
{\bf Proof.} It is enough to show that the 
equality in the theorem holds for all 
$t=t_0$ where (\ref{8nov94(2)}) holds.  
Fix such a $t_0$, and let $A\in \Sigma$.  Since the function 
$t\mapsto T_t\mu(A)$ is bounded on $\R$, it 
follows from the properties of the Poisson kernel that
$$P_y*(T_t\mu)(A)\rightarrow T_t\mu(A)\ \ {\rm for\ almost\ all}\ t\in \R.$$
But by Lemma \ref{secondrelation}, we have
$$P_y*(T_t\mu)(A)\rightarrow T_{t-t_0}
\left(g(t_0)\right)(A)\ \ {\rm for\ all}\ t\in \R.$$
Hence
$$T_t\mu (A)=T_{t-t_0}(g(t_0))(A)$$
for almost all $t\in \R$.  
It is clear from Lemma \ref{secondrelation} that the measure
$T_{-t_0}g(t_0)$ is weakly analytic, since it is the
strong limit in $M(\Sigma)$ of weakly analytic measures.
Applying hypothesis $(A)$, we infer that
$$\mu=T_{-t_0}g(t_0).$$
Applying $T_{t_0}$ to both 
sides of the last equality completes the proof.

\bigskip
From Theorem \ref{minimaintheorem} we can derive several interesting 
properties of analytic measures, which, as the reader may check, 
are equivalent to the F.\ and M.\ Riesz Theorem in the classical  
setting.  We start with a property of the 
Poisson integral of weakly analytic measures. 
\begin{poissonlimit}
Let $T$ and $\mu$ be as in Theorem \ref{minimaintheorem}.
Then,
$$\lim_{y\rightarrow 0} P_y*\mu= \mu$$
in the $M(\Sigma)$-norm.
\label{poissonlimit}
\end{poissonlimit}
{\bf Proof.} 
Let $t_0\in \R$ be such that $P_y*g(t_0)\rightarrow g(t_0)$
in the $M(\Sigma)$-norm (recall (\ref{8nov94(2)})).  
We have
$$T_{-t_0}\left( P_y*g(t_0)\right)\rightarrow T_{-t_0}\left(g(t_0)\right)$$
in the $M(\Sigma)$-norm.  But
$T_{-t_0}\left( P_y*g(t_0)\right)=P_y*\mu$,
and 
$T_{-t_0} g(t_0)=\mu$, and the theorem follows.

\bigskip
The following generalizes Theorem 4 of Forelli \cite{forelli} 
which in turn is a generalization of Theorem (3.1) 
of de Leeuw and Glicksberg \cite{dg}.  
\begin{continuity}
Let $\mu$ and $T$ be as in Theorem \ref{poissonlimit}.  
Then the mapping $t\mapsto T_t\mu$ is uniformly continuous 
from $\R$ into $M(\Sigma)$.
\label{continuity}
\end{continuity}
{\bf Proof.}  It is easily seen that for each $y>0$, the map 
$t \mapsto P_y*g(t)$\ is 
continuous.  By
Lemma \ref{firstrelation}, it follows that the map $t \mapsto P_y*T_t\mu
= T_t(P_y*\mu)$\ is continuous.  By 
Theorem \ref{poissonlimit} and (\ref{uniformlybded}), we see that
$T_t(P_y*\mu) \to T_t\mu$\ uniformly in $t$, and the result follows.
\begin{corbigmaintheorem2}
Let $T$ be a representation of $\R$ acting on 
$M(\Sigma)$ and satisfying hypothesis $(A)$,
(\ref{uniformlybded}), and
(\ref{adjointcond}).  
Suppose that $P$ is 
a bounded linear operator from $M(\Sigma)$ into itself 
that commutes with $T_t$\ for each $t\in \R$.  If $\mu$ 
is a weakly analytic measure in $M(\Sigma)$, then
${\rm spec}_T(P\mu)$\ is contained in ${\rm spec}_T\mu$.  In
particular, $P\mu$\ is weakly analytic.
\label{corbigmaintheorem2}
\end{corbigmaintheorem2}
{\bf Proof.}  Clearly, it is sufficient to show that
${\cal J}(\mu) \subset {\cal J}(P\mu)$.  So, suppose that
$f \in {\cal J}(\mu)$.  Define the measure
$ \nu = f*_T \mu $.  Now, the map $t\mapsto T_t\mu$\ is Bochner 
measurable, and hence the map
$t\mapsto f(t) T_{-t}\mu$\ is Bochner integrable.  By properties of the
Bochner integral, it follows that
$$ \int_\R f(t) T_{-t} \mu dt = \nu .$$
From the definition of ${\cal J}(\mu)$, we see
that for each $A\in\Sigma$
$$ T_s\nu(A) = \int_\R f(t) T_{s-t}\mu(A) dt = 0 \quad\mbox{a.e.\ $s$.} $$
Hence by hypothesis~($A$), it follows that $\nu = 0$.
Thus, once again, using the properties of the Bochner integral, we see
that
\begin{eqnarray*}
\int_\R f(s-t) T_t P\mu dt
&=&
\int_\R f(t) T_{s-t} P\mu dt \\
&=& P T_s \left(\int_\R f(t) T_{-t} \mu dt\right) \\
&=& P T_s \nu = 0 .
\end{eqnarray*}
Hence, for every $A\in\Sigma$, we have that
$$ \int_\R f(s-t) T_t (P\mu)(A) dt = 0,$$
that is, $f \in {\cal J}(P\mu)$.

\begin{8nov94remark}
{\rm 
When $\Sigma$ is countably generated, using 
Theorem \ref{riesz1}, instead of the Main Theorem of
this section, we can derive a version
of Theorem \ref{corbigmaintheorem2}
without the additional condition on the 
adjoint (\ref{adjointcond}), and more interestingly, 
without assuming that $T$ is a representation by a group of 
isomorphisms on $M(\Sigma)$.  The hypotheses of Theorem \ref{riesz1}
are enough to derive these results. 
}
\label{8nov94remark}
\end{8nov94remark}
We end this section with the proofs of the properties of
convolutions that we stated at the outset of this section.
Throughout the rest of this section, we
use the following notation:
$\mu$ is a weakly measurable element in $M(\Sigma)$;
$\nu$ and $\sigma$ are regular Borel measures in $M(\R)$;
$T=(T_t)_{t\in\R}$ is a one-parameter group of
operators on $M(\Sigma)$ satisfying hypothesis $(A)$ and
such that (\ref{uniformlybded}) and (\ref{adjointcond}) hold. 
The convolution of $\mu$ and $\nu$ is defined on the 
sigma-algebra $\Sigma$ by
\begin{equation}
\nu*_T\mu(A)=\int_\R T_{-t}\mu(A)d \nu(t),
\label{convolution}
\end{equation}
for all $A\in \Sigma$.  When there is no risk
of confusion we will simply write $\nu*\mu$ for $\nu*_T\mu$.

Using dominated 
convergence, it is easy to check that (\ref{convolution}) 
defines a measure in $M(\Sigma)$, and that 
$\|\nu*\mu\|\leq c \|\mu\| \|\nu\|$, 
where $c$ is the as in (\ref{uniformlybded}).  
\begin{applem1}
Suppose that $f\in \cL^\infty(\Sigma)$.  
Then the mapping
$t\mapsto \int_\Omega fd(T_t\mu)$
is Lebesgue measurable on $\R$.  Furthermore, 
\begin{equation}
\int_\R\int_\Omega f d(T_{-s})\mu d\nu(s)=\int_\Omega f d\nu*\mu.
\label{8nov94(3)}
\end{equation}
\label{applem1}
\end{applem1}
{\bf Proof.}  It is sufficient to prove the lemma in the case when $f$\
is a simple function, and then it is obvious.

\begin{appcor1}
For all $t\in \R$, we have
$$T_t(\nu*\mu)=\nu*(T_t\mu).$$ 
Moreover, the measure $\nu*\mu$ is weakly measurable.
\label{appcor1}
\end{appcor1}
{\bf Proof.}  For $A\in \Sigma$, we have
\begin{eqnarray*}
\nu *(T_t\mu)(A)        &=&     \int_\R (T_{-s +t}\mu)(A)d\nu(s)\\
			&=&     \int_\R\int_\Omega 
				T^*_t 1_A dT_{-s}\mu d\nu (s)\\
			&=&     \int_\R T^*_t1_A d\nu*\mu \ \ 
				({\rm by\ Lemma\ \ref{applem1}})\\
			&=&     \int_\R 1_A dT_t(\nu*\mu)\\
			&=&     T_t(\nu*\mu)(A).
\end{eqnarray*} 
To prove the second assertion, 
note that 
$t\mapsto T_t(\nu*\mu)(A)=\nu*(s\mapsto T_s\mu(A))(t)$, 
and so the function $t\mapsto T_t(\nu*\mu)(A)$ 
is Lebesgue measurable, being the convolution of a measure 
in $M(\R)$ and a bounded measurable function on $\R$.

\begin{appcor2}
With the above notation, we have
$$(\sigma*\nu)*\mu=\sigma*(\nu*\mu).$$
\label{appcor2}
\end{appcor2}
{\bf Proof.}  For $A\in \Sigma$, we have
\begin{eqnarray*}
(\sigma*\nu)*\mu(A)     &=&     
\int_\R(T_{-s}\mu)(A)d(\sigma*\nu)(s)\\
			&=&     
\int_\R\int_G(T_{(-s-t)}\mu)(A)d\nu(t)d\sigma(s)
	       \ \ ({\rm by\ \cite[Theorem\ (19.10)]{hr1}})\\
			&=&     \int_\R (\nu*(T_{-s}\mu))(A)d\sigma(s)\\
			&=&     \int_\R T_{-s} (\nu*\mu)(A)d\sigma(s)
				\ ({\rm by\ Corollary\  \ref{appcor1}})\\
			&=&     \sigma*(\nu*\mu)(A),
\end{eqnarray*}
and the lemma follows.


\section{Lebesgue decomposition of analytic measures}
\newtheorem{defqi}{Definition}[section]
\newtheorem{decomposition}[defqi]{Theorem}
\newtheorem{8nov94remark2}[defqi]{Remark}
\newtheorem{notstronganalytic}[defqi]{Example}
\newtheorem{fandmriesz2}[defqi]{Theorem}
In their extension of the F.\ and M.\ Riesz Theorem
to compact abelian groups, Helson
and Lowdenslager \cite{hl1} realized that while an analytic
measure may not be absolutely continuous with respect to
Haar measure, its absolutely continuous part and singular part
are both analytic.  This property was then generalized by
de Leeuw and Glicksberg \cite{dg} and Forelli \cite{forelli} 
to the Lebesgue decomposition of analytic measures with respect
to quasi-invariant measures, which take the place of
Haar measures on arbitrary measure spaces.  
In this section, we derive our version of this result as a simple
corollary of Theorem \ref{corbigmaintheorem2}.  We then 
derive a version of this theorem in the case when $\Sigma$ consists
of the Baire subsets of a topological space,
without using the group property of the 
representation.

The setting for this section is as follows. 
Let  
$T=(T_t)_{t\in\R}$ denote a one-parameter group given by 
mappings of a sigma algebra  $\Sigma$.  With a slight abuse of notation,
we will write
$$T_t\mu(A)=\mu(T_t(A))$$
for all $t\in\R$, all $A\in \Sigma$, and all 
$\mu\in M(\Sigma)$.  Note that
conditions (\ref{uniformlybded}) and (\ref{adjointcond}) are satisfied.  
In addition to these properties,
we suppose that $T$ satisfies hypothesis $(A)$, and so the results of the
previous section can be applied in our present setting.
\begin{defqi}
Let $T$ be as above, and let 
$\nu\in M(\Sigma)$ be weakly measurable.  
We say that $\nu$ is quasi-invariant if, 
for all $t\in\R$, $\nu$ 
and $T_t\nu$ are mutually absolutely continuous.
\label{defqi}
\end{defqi}

The following is a generalization to arbitrary measure spaces of 
Theorem 5 of Forelli \cite{forelli}.
\begin{decomposition}
Let $T$ be as above, and
let $\mu$ and $\sigma$ be weakly measurable in $M(\Sigma)$ such that
$\sigma$ is quasi-invariant and $\mu$ is weakly analytic.
Write $\mu=\mu_a+\mu_s$ for the Lebesgue decomposition of 
$\mu$ with respect to $\sigma$.  Then the spectra of 
$\mu_a$ and $\mu_s$ are contained
in the spectrum of $\mu$.  In particular, $\mu_a$ and 
$\mu_s$ are both
weakly analytic in $M(\Sigma)$.
\label{decomposition}
\end{decomposition}
{\bf Proof.}  Define $P$ on $M(\Sigma)$ by
$P(\eta)=\eta_s$, where $\eta_s$
is the singular part of $\eta$ in its Lebesgue 
decomposition with respect to $\sigma$.  

It is easy to see that the quasi-invariance of
$\nu$ is equivalent to the fact that for all
 $A\in \Sigma$, we have $|\nu|(A)=0$ if and
only if $T_t|\nu|(A)=0$ for all $t\in\R$.  Consequently, 
$$P\circ T_t(\eta)=T_t\circ P(\eta).$$
Now
apply Theorem \ref{corbigmaintheorem2}.
\begin{notstronganalytic}
{\rm Consider Examples \ref{ex3} (a) and (b).  
Clearly, the measure $\mu = (\nu -\delta_0)$ 
is weakly analytic with respect to the 
representations $T^\alpha$, for any $\alpha$, 
since, for every $A\in \Sigma$, 
$T_t^{\alpha} \mu (A)=0$ for almost 
all $t\in \R$, and so the
function $t\mapsto T^\alpha_t\mu(A)$
is trivially in $H^\infty(\R)$.  
However, $\mu_s = -\delta_0$, and hence
$T_t^\alpha \mu_s = -e^{i\alpha} \delta_t$, which is not 
weakly analytic if $\alpha < 0$.}
\label{notstronganalytic}
\end{notstronganalytic}
\begin{8nov94remark2}
{\rm Using 
Remark \ref{8nov94remark}, we see that
Theorem \ref{decomposition} also holds under the hypotheses 
of Theorem \ref{riesz1}. 
}
\label{8nov94remark2}
\end{8nov94remark2}
We can use Remark \ref{8nov94remark2} to show that, on topological
spaces where the action of $\R$ is given by 
jointly continuous point mappings of the underlying
space,
Theorem \ref{decomposition} holds even if we dispense with the group
property of the representation.
\begin{fandmriesz2}
Let $\Omega$ be a topological space, and let $(T_t)_{t\in\R}$ be 
a family of homeomorphisms of $\Omega$ such that 
$(t, \omega)\mapsto T_t\omega$ is jointly continuous. 
Suppose that $\mu$ and $\nu$ are Baire measures
such that $\nu$ is quasi-invariant, and write $\mu=\mu_a+\mu_s$
for the Lebesgue decomposition of $\mu$ with respect to $\nu$.
If $\mu$ is weakly analytic, then the spectra of $\mu_a$ and $\mu_s$
are contained in the spectrum of $\mu$.  In particular, $\mu_a$ and 
$\mu_s$ are both
weakly analytic in $M(\Sigma)$.
\label{fandmriesz2}
\end{fandmriesz2}
{\bf Proof.}  It is enough to consider
$\mu_s$.  Let $\Sigma$ denote the sigma algebra of 
Baire subsets of $\Omega$, and let $A\in \Sigma$.
We want to show that the 
mapping $t\mapsto T_t\mu_s(A)$ is in $H^\infty(\R)$.
We will reduce the problem to a countably generated
subsigma algebra of $\Sigma$ that depends
on $A$, then use Remark \ref{8nov94remark2}.

A simple argument shows that
for each $C\in \Sigma$, there exist a countable collection
of continuous function $\{f_n:\Omega\to\R\}$\ such that $C$\ is contained
in the minimal sigma-algebra for which the functions $f_n$\ are measurable.
Furthermore, for each continuous function, we have that 
$f\circ T_t \to f\circ T_{t_0}$\ pointwise as $t\to t_0$.  Hence, we see that
$C$\ is an element of 
\begin{eqnarray*}
\Sigma(C) 
&=&
\sigma \left(
\left( f_n\circ T_{r_{1}}^{\pm 1} 
T_{r_{2}}^{\pm 1}
\ldots
T_{r_{k}}^{\pm 1}\right)^{-1} (-r,\infty):n\in \N,\,r_1,r_2,\dots,r\in\R
\right) \\
&=&
\sigma\left( 
\left(f_n\circ T_{r_{1}}^{\pm 1} 
T_{r_{2}}^{\pm 1}
\ldots
T_{r_{k}}^{\pm 1} \right)^{-1}
(-r,\infty):n\in \N,\,r_1,r_2,\dots,r\in\Q\right) .
\end{eqnarray*}
Clearly $\Sigma(C)$\ is countably generated, and
invariant under $T_t$\ and $T_t^{-1}$\ for all
$t \in \R$.

Let $B$ denote the support of $\mu_s$, and let 
$\Sigma(A,B)=\sigma\left( \Sigma(A),\Sigma(B)\right)$.
Again, we have that $\Sigma(A,B)$ is countably generated
and invariant by all $T_t^{\pm 1}$.  
Let $\left. \mu_s\right|_{\Sigma(A,B)}$, 
$\left. \mu\right|_{\Sigma(A,B)}$,
and
$\left. \nu\right|_{\Sigma(A,B)}$
denote the restrictions of $\mu_s$, $\mu$, and 
$\nu$ to the
$\Sigma(A,B)$, respectively.
It is clear that
$\left. \nu\right|_{\Sigma(A,B)}$
is quasi-invariant and that 
$\left. \mu\right|_{\Sigma(A,B)}$
is weakly analytic, where here we are restricting
the definitions to the smaller sigma algebra $\Sigma(A,B)$.
By Remark \ref{8nov94remark2},
the measure
$\left( \left. \mu\right|_{\Sigma(A,B)}\right)_s$
is weakly analytic.  But,
$\left. \mu_s\right|_{\Sigma(A,B)}=
\left( \left. \mu\right|_{\Sigma(A,B)}\right)_s$,
and hence $t\mapsto T_t\mu_s(A)$ is in $H^\infty (\R)$,
completing the proof of the theorem.

\section{Quasi-invariance of analytic measures} 
\newtheorem{lemma1qi}{Lemma}[section]
\newtheorem{lemma2qi}[lemma1qi]{Lemma}
\newtheorem{lemma3qi}[lemma1qi]{Lemma}
\newtheorem{thqi}[lemma1qi]{Theorem}
\newtheorem{exqi}[lemma1qi]{Example}
In this last section, 
we use some of the machinery that we have developed in the 
previous sections to give a simpler proof of a result
of de Leeuw, Glicksberg, and Forelli, that asserts
that analytic measures are quasi-invariant.  
We will show by an example that unless the action is restricted
to point mappings of the underlying topological space, 
such a result may fail.
Thus, the various results that we obtained in Sections
2 and 3 for more general representations 
of $\R$ cannot be obtained by the
methods of Forelli \cite{forelli} which imply the 
quasi-invariance of analytic measures.


As in the previous sections, we start by 
describing the setting for our work.
Here $(T_t)_{t\in\R}$ 
denotes a group of homeomorphisms of a topological
space $\Omega$, and $\Sigma$ stands for the Baire 
subsets of $\Omega$.  
Given 
a Baire measure $\mu$, we let $T_t\mu$ be the 
measure defined on the Baire subsets $A\in \Sigma$ 
by $(T_t\mu)(A)=\mu(T_t A)$.   
Applying
Theorem \ref{poissonlimit}, we have that $P_y*\mu\rightarrow \mu$ 
in the $M(\Sigma)$ norm.  Using the 
Jordan decomposition of the measure $\mu$, 
we see that $P_y*|\mu |\rightarrow |\mu |$.    
Hence, from the proof of Theorem \ref{continuity}, we see 
that the mapping $t\mapsto T_t|\mu |$ is also continuous.  
Define a measure $\nu$ in $M(\Sigma)$, by
\begin{equation}
\nu(A)=\frac{1}{\pi} \int^\infty_{-\infty} T_t|\mu| 
(A) \frac{dt}{1+t^2},
\label{poissonintegral}
\end{equation}
for all $A\in \Sigma$.  Note that $\nu=P_1 *|\mu|$.
\begin{lemma1qi}
For all $t\in \R$, we have $T_t|\mu|<<\nu$, and hence $T_t\mu<<\nu$.
\label{lemma1qi}
\end{lemma1qi}
{\bf Proof.}  Let $A\in \Sigma$.  Since the mapping 
$t\mapsto T_t|\mu |(A)$ is continuous and nonnegative, the 
lemma follows easily.

\begin{lemma2qi}
Let $h(t,\omega)$ denote the Radon-Nikod\'ym 
derivative of $T_t\mu$ with respect to $\nu$.  Then 
the mapping $t\mapsto h(t,\cdot)$ is 
continuous from $\R$ into $L^1(\nu)$.  
Consequently it is Bochner measurable, and 
hence $(t,\omega)\mapsto h(t,\omega)$ is jointly 
measurable on $\R\times\Omega$.
\label{lemma2qi}
\end{lemma2qi}
{\bf Proof.}  The first part of the lemma follows once 
we establish that the real and imaginary parts, 
$t\mapsto \Re h(t,\cdot)$ 
and $t\mapsto \Im h(t,\cdot)$, are continuous.  We deal 
with the first function only; the second is 
handled similarly.  Let 
$B=\{\omega\in \Omega:\ \ \Re ( h(t,\omega)-h(t',\omega))>0\}$. 
We have 
\begin{eqnarray*}
\lefteqn{
\|\Re (h(t,\cdot)-h(t',\cdot))\|_{L^1(\nu)}     =     
	\int_\Omega  | \Re (h(t,\omega)-h(t',\omega))|d\nu(\omega)}\\
					&=&
 \int_B \Re (h(t,\omega)-h(t',\omega))d\nu(\omega)+\int_{\Omega\setminus B} 
	\Re (h(t',\omega)-h(t,\omega))d\nu(\omega)\\
					&\leq&
	|T_t(|\mu |) (B)-T_{t'}(|\mu |) (B)|+ 
	|T_t (|\mu |) (\Omega\setminus B)-T_{t'}(|\mu |)
	(\Omega\setminus B)|.
\end{eqnarray*}
The continuity follows now from Lemma \ref{lemma1qi}.  
To complete the proof of the lemma, note that since 
$t\mapsto h(t,\cdot)$ is Bochner measurable, it 
is the limit of simple functions $h_n(t,\cdot)$ 
each of which is jointly measurable on $\R\times \Omega$.  

\bigskip
One more property of the function $h$ is needed before 
establishing the quasi-invariance of analytic measures.
\begin{lemma3qi}
Suppose that $\mu\in M(\Sigma)$ 
is weakly analytic, and let $h(t,\omega)$ 
be as in Lemma \ref{lemma2qi}.  Then there is 
a Baire subset $\Omega_0$ of $\Omega$ such 
that $\nu (\Omega\setminus\Omega_0)=0$, and 
for all $\omega\in \Omega_0$, the 
function $t\mapsto h(t,\omega)$ is in $H^\infty(\R)$. 
\label{lemma3qi}
\end{lemma3qi}
{\bf Proof.} By
Lemma \ref{lemma2qi}, the function $t\mapsto h(t,\cdot)$
is bounded and Bochner measurable.  Hence, for any $g(t)$ in $H^1(\R)$, the  
function $h(t,\omega) g(t)$ is in $L^1(\R,L^1(\nu))$.  
Moreover, for any Baire subset $A\in \Sigma$, we have
\begin{eqnarray}
\int_A \int_\R h(t,\omega) g(t) dt d\nu(\omega) &=&
		\int_\R \int_A h(t,\omega) g(t) d\nu(\omega)dt \nonumber\\
		&=& 
		\int_\R T_t\mu(A) g(t) dt=0
\label{lemm3qi1}
\end{eqnarray}
since $\mu$ is weakly analytic.  Since (\ref{lemm3qi1}) holds 
for every $A\in \Sigma$, we conclude that the function 
$\omega \mapsto \int_\R h(t,\omega)g(t)dt=0$ for $\nu$-almost 
all $\omega$.  Hence there is a subset $\Omega_g\in \Sigma$ 
such that $\nu (\Omega\setminus\Omega_g)=0$ 
and $\int_\R h(t,\omega)g(t)dt=0$ for all 
$\omega\in \Omega_g$.  Since $H^1(\R)$ is 
separable, it contains a countable dense 
subset, say $\{g_n\}$.  Let 
$ \Omega_0=\bigcap \Omega_{g_n}$.  
Then, $\nu(\Omega\setminus\Omega_0)=0$, 
and for $\omega\in \Omega_0$ and all 
$g\in H^1(\R)$ we have that $\int_\R h(t,\omega)g(t)dt=0$ 
which proves the lemma.

\bigskip

We now come to the main result of this section.  In addition to the 
preliminary lemmas that we have just established, the proof 
uses the fact that a function in $H^\infty(\R)$ is either zero almost
everywhere or is not zero almost everywhere.
\begin{thqi}
Suppose that $t\mapsto T_t$ is a one-parameter 
group of homeomorphisms of a topological space $\Omega$ 
with the property that $(t,\omega)\mapsto T_t\omega$ 
is continuous, and let $\Sigma$ denote the sigma-algebra of 
Baire subsets of $\Omega$.  Suppose that $\mu\in M(\Sigma)$ 
is weakly analytic.  Then $\mu$ is quasi-invariant.
\label{thqi}
\end{thqi}
{\bf Proof.}  We use the notation of the previous lemma.  
For $\omega\in\Omega_0$ and $A\in\Sigma$, we have
\begin{eqnarray}
T_t|\mu|(A)     &=&     T_{t+s}|\mu| \left( T_{-s}A\right)
			\nonumber \\
		&=&     \int_{T_{-s} A} |h(t+s,\omega)|d\nu(\omega) .
			\nonumber 
\end{eqnarray}
Hence,
\begin{eqnarray}
T_t|\mu|(A)     &=&     \int_\R T_t|\mu|(A) \frac{ds}{\pi(1+s^2)}
			\nonumber \\
		&=&     
\int_\R \int_{T_{-s} A} |h(t+s,\omega)|d\nu(\omega)\frac{ds}{\pi(1+s^2)}
			\nonumber \\
		&=&     
\int_\Omega \int_{C_\omega}|h(t+s,\omega)|
\frac{ds}{\pi(1+s^2)}d\nu(\omega)
			\nonumber \\
		&=&     
\int_{\Omega_0} \int_{C_\omega}|h(t+s,\omega)|
\frac{ds}{\pi(1+s^2)}d\nu(\omega),                
\label{thqi1}
\end{eqnarray}
where $C_\omega=\{s\in\R:\ \ \omega\in T_{-s}A\}$.  
Since for $\omega\in\Omega_0$, the function 
$t\mapsto h(t,\omega)$ is in $H^\infty(\R)$, it 
follows that this function is either zero $t$-a.e., or 
not zero $t$-a.e.  Let 
$\Omega_1=\{\omega\in\Omega:\ \ h(t,\omega)=0, 
{\rm for\ almost\ all}\ t\in\R\}$.  
Then 
$\Omega_1\in\Sigma$, and from (\ref{thqi1}) we have that
\begin{equation}
T_t|\mu|(A)=\int_{\Omega_0\setminus\Omega_1} 
\int_{C_\omega}|h(t+s,\omega)|\frac{ds}{\pi(1+s^2)}d\nu(\omega).
\label{thqi2}
\end{equation}
Hence, $T_t|\mu|(A)=0$ if and only if 
$$\int_{C_\omega}|h(t+s,\omega)|\frac{ds}{\pi(1+s^2)}=0$$
for $\nu$-almost all $\omega\in \Omega_0\setminus\Omega_1$.  
Since the integrand is strictly positive except on a set of zero measure, 
this happens 
if and only if the Lebesgue measure of $C_\omega$ is zero for 
$\nu$-almost all $\omega\in \Omega_0\setminus\Omega_1$.  
But since this last condition does not depend on $t$, we see 
that $T_t|\mu|(A)=0$ if and only if $|\mu|(A)=0$.

\bigskip

That Theorem \ref{thqi} does not hold for more 
general representations is demonstrated by the following example.
\begin{exqi}
{\rm
Let $\Omega=\{0,1\}$, and let $\Sigma$ 
consist of the power set of $\Omega$.  Denote by 
$\delta_0$ and $\delta_1$ the point masses at 0 
and 1, respectively.  Then $ \delta_0$ and $\delta_1$ 
form a basis for $M(\Sigma)$, and every element 
in $M(\Sigma)$ will be represented as a vector 
in this basis.  For $t\in\R$, define $T_t$ by the matrix
\[ T_t=\left(
\begin{array}{cc}
e^{4 i t}\cos t & e^{4 i t}\sin t \\
- e^{4 i t}\sin t & e^{4 i t}\cos t 
\end{array}
\right).\]
Note that $T$ satisfies hypothesis $(A)$.  
Also, it is easy to verify that $\delta_0$ is 
weakly analytic.  However, $T_{\pi/2}\delta_0=-\delta_1$, 
and $\delta_0$ and $\delta_1$ are not mutually absolutely continuous.
Hence $\delta_0$ is weakly analytic but not quasi-invariant.
}
\label{exqi}
\end{exqi}
{\bf Final Remarks}  The approach that we took to the 
F.\ and M.\ Riesz Theorem can be carried out in the more general setting
where $\R$ is replaced by any locally compact abelian group $G$
with an ordered dual group $\Gamma$, and where the notion
of analyticity is defined as in \cite{hl1} using the order
structure on the dual group.
With the exceptions of $G=\R$ and $G=\T$, the main result of
\S 3, concerning the Bochner measurability of $t\mapsto T_t\mu$, fails
even in the nice setting of a regular action of $G$ by translation in 
$M(G)$.  However, we can prove a weaker result that states that 
a weakly analytic measure is strongly analytic, that is, whenever
$\delta \in M(\Sigma)^*$ (the Banach dual of $M(\Sigma)$), then the map $t\mapsto\delta(T_t \mu)$\ is
in $H^\infty(G)$.  This result, in turn, implies the
versions of the F.\ and M.\ Riesz 
Theorems proved by Helson and Lowdenslager \cite{hl1} 
for actions of locally compact abelian groups
on abstract measure spaces.  This work 
is done in a separate paper by the authors.

\bigskip\noindent
{\bf Acknowledgements}  The work of the authors was supported
by separate grants from the National Science Foundation (U. S. A. ).

\end{document}